\documentclass[bjps]{imsart}
\pdfoutput=1

\RequirePackage[OT1]{fontenc}
\usepackage{amsthm,amsmath,natbib}

\usepackage{amssymb}
\RequirePackage{hyperref}

\startlocaldefs
\newcommand{\BS}{\mathbb S}
\newcommand{\RR}{\mathbb R}
\newcommand{\SOd}{\mathit{SO(d)}}
\newcommand{\Od}{\mathit{O(d)}}
\newcommand{\SOtwo}{\mathop{SO(2)}}
\newcommand{\Otwo}{\mathop{O(2)}}
\newcommand{\SOthree}{\mathop{SO(3)}}
\newcommand{\Othree}{\mathop{O(3)}}

\newcommand{\Ga}{\Gamma}
\newcommand{\Th}{\Theta}

\renewcommand{\d}{{\mathrm{d}}}
\newcommand{\ka}{\kappa}
\newcommand{\si}{\sigma}
\renewcommand{\th}{\theta}

\newcommand{\CAY}{\mathrm{\,CAY}}

\newcommand{\E}{\mathrm E}

\endlocaldefs

\startlocaldefs
\numberwithin{equation}{section}
\theoremstyle{plain}
\newtheorem{defn}{Definition}
\newtheorem{prop}{Proposition}
\newtheorem{thm}{Theorem}
\newtheorem{rem}{Remark}
\endlocaldefs

\begin{document}

\begin{frontmatter}
\title{Fake Uniformity in a Shape Inversion Formula}
\runtitle{Fake Uniformity in Shape Inversion}

\begin{aug}
\author{\fnms{Christian} \snm{Rau}\thanksref{a}\ead[label=e1]{christianrau080@gmail.com}}


\runauthor{C. Rau}

\affiliation[a]{Department of Mathematics, Shantou University}

\address{Christian Rau, Department of Mathematics,\\
Shantou University, Shantou Guangdong 515063, P.R.~China \\
\printead{e1}}

\end{aug}

\begin{abstract}
We revisit a shape inversion formula derived by Panaretos in the context of a particle density estimation problem
with unknown rotation of the particle. A distribution is presented which imitates, or `fakes', the uniformity or
Haar distribution that is part of that formula.  
\end{abstract}




\end{frontmatter}

\section{Introduction} \label{sec:introduction}
Stochastic geometry makes extensive use of uniform, or Haar distributed,
rotations; for example, in constructing random geometric objects (hyperplanes, polytopes, 
or other objects) that have a uniform orientation. The uniformity often makes it
possible to reconstruct a three-dimensional object from two-dimensional projections or
sections; such tasks belong to the realm of stereology, an area which is connected to both
stochastic geometry and spatial statistics.
This paper is a cautionary note, illustrating that the said uniformity 
may, in the case of a natural geometric functional, be imitated, or `faked', 
 by a non-uniform distribution.

The natural functional that we consider here is a central ingredient in a
 {\sl shape inversion} formula from~\cite{Pan09}, which he used 
  to tackle a particle reconstruction problem arising from electron microscopy data.
This formula, which is given in~Proposition~\ref{shape-inv-panaretos} below,
allows to recover salient features, which may be referred to as `landmarks',
  of an object in~$\RR^d$ ($d\ge 2$)
from its projections on an arbitrary fixed $(d-1)$-dimensional subspace.
An important complicating feature of the problem is that the 
object is subject to a prior random and unknown uniform rotation before it is imaged. 
 For practically relevant imaging tasks, where the object may be a protein fragment, we have $d=3$, and the subspace
is the imaging plane of the microscope. In this introduction we
 limit ourselves to those notions that are necessary to formulate the 
 shape inversion formula and our theorem, and refer to  
Section~\ref{ramifications} for a precise statement of the statistical model.

The precise definition of uniform or Haar (probability) measure of a rotation is as follows. 
 Recall that  $\SOd$, the group of $d\times d$ rotation matrices, consists of those matrices $A$ satisfying
$A^{\sf T}A=I$, the $d\times d$ identity matrix,  and~$\det A=1$.
A random rotation~$A$ has the Haar distribution, written as $A\sim \mu$, 
 if the distribution of $QA$ is the same as that of $A$,
for any nonrandom~$Q\in\SOd$.

Here we examine how---and indeed, if---the formula given in 
Proposition~\ref{shape-inv-panaretos} changes if
 one replaces~$\mu$ by a less symmetric rotation.
A motivation for this question comes from the fact that 
the group $\SOd$ {\sl acts} on the $d$-dimensional unit sphere $\BS^{d-1} =\{x\in\RR^d:\, \|x\|=1\}$ via matrix-vector multiplication~$x\mapsto Ax$.
In measure theory, this connection between $\SOd$ and $\BS^{d-1}$ is frequently exploited, and the present paper may be seen
as another modest instance. The following fact, which is relevant to this paper, illustrates this connection: 
if $A\sim\mu$ and $v$ is any
nonrandom point on~$\BS^{d-1}$, then $Av$ is has the uniform distribution $\si$ (normalized Lebesgue surface area measure)
on~$\BS^{d-1}$. We note in passing that one may conversely construct $\mu$ from~$\si$~\cite[pp.~584--585]{SW08}.

In Section~\ref{fake}, after some preliminaries, we state the shape inversion formula
of~\cite{Pan09}, and formulate  the fake uniformity problem. In Section~\ref{main-result}, the 
Cayley distribution with $\ka=1$ is found to be a distribution which fakes uniformity for the Gram matrix functional. Finally, 
Section~\ref{ramifications} gives the statistical model in detail, and discusses the impact of 
fake uniformity.

\section{(Closely) faking the value of a functional} \label{fake}
We begin with describing a cousin  of the problem studied in the present note. 
This cousin problem is mathematically more sophisticated but, on the other hand, does not make any reference to rotations.
For a compact convex set $K\subset \RR^d$ and $u\in\BS^{d-1}$, 
denote by $K_u$  the orthogonal projection of $K$ onto the (hyper)plane 
$\{x\in\RR^d:\, \langle x,u\rangle=0\}$, where $\langle \cdot,\cdot\rangle$ is the usual inner product. 
Let $v_{d-1}(\cdot)$ denote $(d-1)$-dimensional volume.
 Consider the functional 
\begin{equation} \label{cauchy-general}
s(\Phi,K) =  \int_{\BS^{d-1}} v_{d-1}(K_u)\,\Phi(u)\,\d \si(u)\,,
\end{equation}
where $\Phi$ is an integrable function on~$\BS^{d-1}$. Since we consider projections, we may assume $\Phi$ 
to satisfy $\Phi(u)=\Phi(-u)$ for all $u$. 
For the constant function $\Phi\equiv 1$, Cauchy's surface area formula~\cite[p.~45]{Gro96} yields that 
$s(\Phi,K)$ is the Lebesgue surface area of~$K$. (Note that our $\si$, unlike in~\cite{Gro96}, is already normalised.)
 In~\cite[pp.~297ff]{Gro96}, the following {\sl inverse problem\/} was considered: 
If $s(\Phi,K)$ is close to $s(1,K)$ for all~$K$ with surface area bounded by some constant, can it be inferred that $\Phi$ is
close to~$1$? The answer is negative, unless smoothness assumptions on $\Phi$ are imposed;~\cite{Gro96} proved
this with tools from harmonic analysis on~$\BS^{d-1}$. 

The `fake uniformity' from the title of the present note is a negative answer to the cousin problem on $\SOd$ which we shall formulate
at the end of this section, and 
a manifestation of a comment on ill-definedness in~\cite[p.~3303]{Pan09}. Our result is not of a limiting nature, and will not require
 harmonic analysis for its proof.

The shape inversion formula of~\cite{Pan09} is as follows. 
Let $V$ be any real $d\times \ell$ matrix  ($d\ge2,\ell\ge 1$), 
and~$H=\mathrm{diag}(1,\ldots,1,0)$ be the projection matrix (with respect to the standard basis of~$\RR^d$)
which zeroes the last component of a vector in~$\RR^d$. We interpret the columns of $V$ as the location vectors of the 
so-called \textit{landmarks} associated with the unknown particle to be reconstructed. 
\begin{prop}[Panaretos, 2009, Theorem~4.1] \label{shape-inv-panaretos}
\begin{equation*} 
\int_{\SOd} \mathrm{Gram}(HAV)\,\mu(\d A) = \frac{d-1}{d}\,\mathrm{Gram}(V)\,,
\end{equation*}
where $\mathrm{Gram}(W)=W^{\sf T} W$ for any matrix~$W$, so that the entry with index $(i,j)$ in
$\mathrm{Gram}(W)$ is the inner product
of the $i$'th and $j$'th columns of $W$. 
 \end{prop}
 Note that $\mathrm{Gram}(AW)=\mathrm{Gram}(W)$ for any $A\in\Od$, the group of 
 $d\times d$ orthogonal matrices. Hence $\mathrm{Gram}(V)$ `nearly'
encodes   \textit{shape} if the latter were understood to be the information 
 which remains if ``we are not interested
in location, orientation or scale of the 
resulting configuration"~\cite[p.~428]{Ken77}---however, information on reflections is lost, and this will be seen
to be the reason for fake uniformity.

Proposition~\ref{shape-inv-panaretos} says that the original shape 
 can be reconstructed from the projected shape. As noted in~\cite[p.~3286]{Pan09}, this feature is shared
 with~\eqref{cauchy-general}. 
 We shall replace $\mu$ in Proposition~\ref{shape-inv-panaretos} by a distribution which has the weaker 
 symmetry property of {\sl conjugation-invariance}. 
 \begin{defn} \label{conj-inv-def}
 A random rotation $A$ is {\rm conjugation-invariant}
if $Q^{\sf T}AQ$ has the same distribution as~$A$, for any nonrandom~$Q\in \SOd$.
\end{defn}
 If we assume that $A$ has
a density~$f$ with respect to $\mu$, then conjugation-invariance may be expressed by 
the requirement that
$f(Q^{\sf T}PQ)=f(P)$ for all~$Q$ as in Definition~\ref{conj-inv-def}. 
We note that conjugation-invariant functions are also called {\sl central}~\cite[p.~132]{Far08}.

In the case $d=3$, conjugation-invariance has the following geometric meaning. Recall that Euler's theorem says that
every rotation in $\RR^3$ has an axis and angle; if the rotation is assumed to be counter-clockwise and in the interval $(0,\pi)$,
then the orientation of the axis is given by a well-defined vector in~$\BS^2$. 
It can be shown that the conjugation-invariant rotations in~$\RR^3$ which have a density with respect to Haar measure $\mu$
are precisely those for which (i)~the oriented rotation axis $U$ is uniformly distributed on~$\BS^2$, 
and (ii)~the rotation angle~$\Th$ has a density, and~$\Th$ and~$U$ are independent~\cite[Thm.~2.2, p.~109]{Schi97}. 

Let us call the Gram matrix $\mathrm{Gram}(W)$
a `functional' of a given matrix~$W$, even though this is an abuse of terminology, since functionals are usually scalar-valued.
(This reservation is not serious, see Remark~\ref{scalar-reduct} below.) 
For the Gram functional, can the role of $\mu$ in Proposition~\ref{shape-inv-panaretos} be faked by another distribution? We shall affirm this 
by stating an offending alternative random rotation in the next section.

\section{Cayley distribution and main theorem} \label{main-result}
The Cayley distribution was introduced in~\cite{Schae97} under the name of de la Vall{\'e} Poussin
distribution, and independently in~\cite{LMR06}; we adopt the name from the latter reference,
as it has also been used in the R package documentation of~\cite{SHG16}. We state the Cayley density for the 
case $d=3$. This is the case where we
can give an explanation through geometry (see the discussion around~\eqref{doubling});  
however, the proof itself, being based on integration
on~$\BS^{d-1}$, carries over to
the general case~$d\ge 3$.  The Cayley distribution is given by the density ($\Ga(\cdot)$ 
is the Gamma function and $\mathrm{tr(\cdot)}$ is the trace) 
\begin{align*}
f^\CAY_\ka(R) &=   \frac{\sqrt{\pi}\,\Ga(\ka+2)}{2^{2\ka}\Ga(\ka+1/2)}(1+\mathrm{tr}(R))^\ka \nonumber \\
&= \frac{\sqrt{\pi}\,\Ga(\ka+2)}{2^{\ka}\Ga(\ka+1/2)}(1+\cos \th)^\ka\,,  \qquad \quad  0 \le \th \le \pi\,,
\end{align*}
which only depends on the rotation angle~$\th$, and thus is conjugation-invariant; 
the parameter $\ka\ge 0$ measures spread around the median $I$, with the case $\ka=0$ 
corresponding to~$\mu$. 
The density of~$\Th$ is~\cite[p.~424]{LMR06}  
\begin{equation}  \label{Cayley-angle}
f^\CAY_{\Th}(\th) = \frac{\Ga(\ka+2)}{\sqrt{\pi}2^\ka\Ga(\ka+1/2)}(1+\cos \th)^\ka(1-\cos \th)\,,  \qquad 0 \le \th \le \pi\,,
\end{equation}
where the factor $(1-\cos \th)$ shows the preference of $\mu$ for large rotations, see also~\cite[Remark 2.4, pp.~102--103]{Schi97}. 
Write $e_k$ for the $k$'th column of $I$ ($k=1,\ldots,d$), where the dimension $d$ will always be clear from the context, 
and write $f_{v_R}$ for the density of $v_R=Re_d$ with respect to $\si$,
where $R$ is any conjugation-invariant distribution with a density with respect to~$\mu$. 
For the Cayley distribution, the density $f_{v_R}$ is known for any~$d$, see~\cite[Prop.~3.3, p.~419]{LMR06};
in particular
\begin{equation} \label{Cayley-last-column}
f^\CAY_{v_R}(w) = 2^{-\ka}\,(\ka+1)\,(1+w_3)^\ka\,, \qquad w=(w_2,w_2,w_3) \in \BS^2\,.
\end{equation}
A density $f$ on $\BS^{d-1}$ that depends only on $w_d$ is called zonal (with respect to $e_d$). 
Then $f(w)={\overline f}(t)$ for a function ${\overline f}(t)$ which is defined on the 
interval~$[-1,1]$, and integrates to a constant (depending on $d$) which can be evaluated
with formula~\eqref{faraut-formula} below.

We now extend Proposition~\ref{shape-inv-panaretos} to conjugation-invariant rotations.

\begin{thm} \label{main-thm}
Let $P\in\SOd$ be a rotation with Haar density $f_P$, and 
such that there exists a nonrandom $M\in\SOd$ such that  $PM^{\sf T}=R$ is conjugation-invariant. Then, with $V$ as 
in Proposition~\ref{shape-inv-panaretos}, 
\begin{equation*}
\int_{\SOd} \mathrm{Gram}(HAV)f_P(A)\mu(\d A) = V^{\sf T}M^{\sf T}(I-D^2)MV\,,
\end{equation*}
where $D^2=\mathrm{diag}\left((1-\tau_2)/(d-1),\ldots,(1-\tau_2)/(d-1),\tau_2\right)$, and 
\begin{equation} \label{tausq-formula}
\tau_2 =  \frac{\Ga\left(\frac{d}{2}\right)}{\sqrt{\pi}\,\Ga\left(\frac{d-1}{2}\right)} 
\int_{-1}^1 t^2(1-t^2)^{(d-3)/2} {\overline f}_{v_{R^{\sf T}}}(t)\,\d t\,. 
\end{equation}
\end{thm}
\begin{proof}
Similar to~\cite[pp.~3285--3286]{Pan09}, one obtains
\begin{align}
\int_{\SOd} \mathrm{Gram}(HAV)f_P(A)\mu(\d A) 
  &= V^{\sf T}M^{\sf T}\cdot \left\{I-\E\left(v_{R^{\sf T}} v_{R^{\sf T}}^{\sf T} \right)\right\}\cdot MV\,,   \label{sec-mom-1}  \\
 \E\left(v_{R^{\sf T}} v_{R^{\sf T}}^{\sf T} \right) &= \left(\int_{\BS^{d-1}}  w_iw_j f_{v_{R^{\sf T}}}(w)\,\si(\d w) \right)_{i,j}\,. \label{sec-mom-2}
\end{align}
That~\eqref{sec-mom-2} vanishes for $i\neq j$
 follows from symmetry considerations applied in conjunction with the 
 following standard integration formula~\cite[Prop.~9.1.2, p.~189]{Far08}, where $\si_0$ is the uniform distribution
on the `equator'~$\BS^{d-1}_0 = \BS^{d-1}\,\cap\,\{x:\, x_d=0\}$ and $g$ an integrable function: 
\begin{align}
\int_{\BS^{d-1}} & g(x)\,\si(\d x)   \label{faraut-formula} \\
&= \frac{\Ga\left(\frac{d}{2}\right)}{\sqrt{\pi}\,\Ga\left(\frac{d-1}{2}\right)}
\int_0^\pi \left( \int_{\BS^{d-1}_0} g\big((\sin \th) u + (\cos \th) e_d\big)\,\si_0(\d u)\right) \sin^{d-1}\th \,\d \th\,. \nonumber 
\end{align}
That same formula applied to the case $i=j=d$, for which $g$ is zonal, yields~\eqref{tausq-formula}.
Symmetry considerations also imply that all entries with $i=j<d$ coincide. Finally, 
 $\E\{\text{tr}(D^2)\}=\E\big(v_{R^{\sf T}} v_{R^{\sf T}}^{\sf T}\big)=1$, hence $ \E\left(v_{R^{\sf T}} v_{R^{\sf T}}^{\sf T} \right)=D^2$, with $D^2$ as in the theorem.
\end{proof}

\begin{rem}\label{law-of-transpose}
In the case $d=3$, the distributions of $R$ and $R^{\sf T}=R^{-1}$ coincide, as follows readily from the 
axis-angle representation of~$\SOthree$ in Section~\ref{fake}, together with the observation that 
a rotation by the amount $\th\in\RR$ around the oriented axis $u$ is the same as a rotation by the amount~$-\th$ 
around~$-u$. However, the distributions of $R$ and $R^{\sf T}$ 
do not coincide in general for $d>3$; see the characterisations in~\cite[Prop.~2, p.~2768]{SLLBM10}. 
\end{rem}

For the Cayley distribution with $d=3$, Theorem~\ref{main-thm}, Remark~\ref{law-of-transpose} 
and~\eqref{Cayley-last-column} give via integration by parts
(with $\mathrm{diag}(a_1,a_2,a_3)=\sum_{i=1}^3 a_i e_i e_i^{\sf T}$)
 \begin{align*}
\int_{\SOthree} \mathrm{Gram}&\{HAV\}f^\CAY_\ka(A) \,\mu(\d A) = \mathrm{Gram}(V) -  \mathrm{Gram}(D^{\text{CAY}}_\ka MV)\,, \nonumber \\
D^{\text{CAY}}_\ka &= \mathrm{diag} \left(\sqrt{\frac{2(\ka+1)}{6+5\ka+\ka^2}},\sqrt{\frac{2(\ka+1)}{6+5\ka+\ka^2}}, \sqrt{\frac{2+\ka+\ka^2}{6+5\ka+\ka^2}}\right)\,. 
\end{align*}
We see that $D^{\text{CAY}}_\ka$ is a scalar matrix
(i.e., with all diagonal entries identical) not only for the case $\ka=0$, but also for~$\ka=1$. 
This latter case is what we call {\sl fake uniformity} of the Gram matrix functional. We may understand its genesis geometrically 
as follows: the nonconstant factor of~\eqref{Cayley-angle} for $\ka=1$ is
\begin{equation} \label{doubling}
(1+\cos \th)(1-\cos\th) = 1-\cos^2\th = \frac{1}{2}(1-\cos 2\th)\,.
\end{equation}
From the formula for the transformation of the density of a real random variable $X$ under scaling $X\mapsto cX$ with $c\in\RR$, 
we conclude from~\eqref{doubling} that the density of $\Th$ with respect to  
Lebesgue measure on $[0,\pi]$ 
in the fake case ($\ka=1$) is obtained by halving~$\Th$ in the Haar case ($\ka=0$). 
Associate each point $u=(u_1,u_2,u_3)$ on the  hemisphere $\{u\in\BS^2:\, u_3\ge 0\}$ 
with its reflection~$T(u) =(u_1,u_2,-u_3)$. Identifying $T$ with its matrix with respect to the canonical basis~$\{e_1,e_2,e_3\}$, we 
observe that $T\in\Othree\setminus\SOthree$, and that~$HTV=HV$ for any landmark matrix~$V$.
Ignoring the set of points on the equator $u_3=0$,
which has measure zero and is negligible, we may now, by suitable choice of either the upper or lower hemisphere, 
produce in the fake case $\ka=1$ the same projected 
configuration as for the Haar case~$\ka=0$. Similar to Remark~\ref{law-of-transpose},
the counter-clockwise angle $\alpha$ changes through reflection 
to $\pi-\alpha$. Also, note that the foregoing description effectively defines a coupling between the cases $\ka=0$ and $\ka=1$.
Alternatively, one may base the argument on the distribution of~$v_R$ from~\eqref{Cayley-last-column}, 
rather than the distribution of $\Th$; this reasoning is, however, not as appealing geometrically.

\begin{rem} \label{scalar-reduct} An examination of the proof of the theorem reveals that 
it suffices to consider the case $\ell=1$$:$ a single landmark vector is enough. This surprising fact, which seems to be at odds with
the notion of `shape' as a configuration of several points, is true because conjugation-invariance is quite a strong symmetry property.
\end{rem}

\begin{rem} While we suspect that there are yet earlier references, Proposition~\ref{shape-inv-panaretos} 
is a straightforward consequence of~\cite[Lemma~2.5, p. 796]{GR04},  in conjunction 
with Remark~\ref{scalar-reduct}. The fake case, however, does not seem to have such a near-precedent.
\end{rem}

\section{Reconstruction from orthogonal views and fake uniformity} \label{ramifications}
In this section we give some insights into the consequences of fake uniformity with regard to
 the tomographic reconstruction problem introduced in~\cite{Pan09} 
  mentioned in Section~\ref{sec:introduction}. 
First we state the random tomography model that he introduced, 
 and recall some issues already known from the case of Haar distributed
rotations. We again limit ourselves to the case~$d=3$.
 
In the random tomography model, 
the unknown particle is construed as a three-dimensional compactly supported probability 
density~$\rho(x)=\rho(x_1,x_2,x_3)$ on~$\RR^3$. 
An observed image of $\rho$ is a (discretized to a regular grid in practice) 
projection of $\rho$ at a random angle, that is, it is given by the compactly supported random field 
\begin{equation} \label{Pan-model}
{\check \rho}(x_1,x_2) = \int_{-\infty}^{+\infty} \rho(U^{-1}x)\,\d x_3\,, 
\end{equation}
where $U\sim \mu$ is called the orientation of~$\rho$; below we write $U\rho(x) = \rho(U^{-1}x)$. 
The stochastic Radon transform of length $N\ge 1$
is a sample $\{{\check \rho}_1,\ldots,{\check \rho}_N\}$
 of $N$ independent and identically distributed (i.i.d.) copies of~${\check \rho}$, generated 
by using a sample $\{U_1,\ldots,U_N\}$ of $N$ i.i.d.~copies of~$U$.

As shown in~\cite{Pan09}, the level of ill-posedness inherent in the problem does not allow the
recovery of $\rho$ itself. However, recovery of $[\rho] = \{A\rho: A\in\Othree\}$, the equivalence class
of $\rho$ which identifies $\rho$ with any rotated or reflected version $A\rho$, {\sl is\/} possible.

In the context of statistical estimation, a low-dimensional parametrization of the particle $\rho$ and 
an arbitrary projection~$\check \rho$ is essential. To this end, 
the roughly spherical `blobs' that are evident in typical images of protein fragments
led~\cite{Pan09} to approximate~$\rho$ by a finite Gaussian mixture with fixed isotropic 
covariance matrices. 
Each component mean is termed a landmark; 
  by themselves, the landmarks 
   give a rough but useful approximation of~$\rho$ or~${\check \rho}$. 
The landmarks are encoded in the columns of the matrix $V$ in Proposition~\ref{shape-inv-panaretos} and Theorem~\ref{main-thm}.
The reason why one focuses on~$[\rho]$ rather than $\rho$ in the reconstruction is also the reason
why one aims to 
reconstruct~$\text {Gram}(V)$ rather than~$V$. 
Proposition~\ref{shape-inv-panaretos} and Theorem~\ref{main-thm} provide the connection between
the Gram matrices for the original three-dimensional and the projected two-dimensional landmarks.

In the context of devising statistical procedures from~Proposition~\ref{shape-inv-panaretos}
in conjunction with (suitable versions of) the law of large numbers and 
central  limit theorem, two issues stand out, as
stated in~\cite[Sec.s~5.1--5.2, pp.~2586--2589]{PK11}. We give them in reverse order.

The second issue was already noted in~\cite[p.~3286]{Pan09} and is
discussed further in~\cite{PK11}: for averaging Gram matrices, one needs to associate
individual landmarks across different images, even though individual labels are {\sl a priori}
not encoded in the averaged Gram matrix.  The randomness in the model is not related to 
this point (except through issues that may arise from the possibly different number of images
at different `views', where `view' is defined in the next paragraph), and hence neither is fake
uniformity.

 The first issue from~\cite{PK11} is whether one may use fewer projections in the
 law-of-large-numbers approximation of the integral on the left-hand side of the equation
 in~Proposition~\ref{shape-inv-panaretos}.
We repeatedly make use, as done in~\cite{PK11}, of the notion of `view' of
the particle. This should be carefully distinguished from the notion of orientation defined around~\eqref{Pan-model}. 
 The notion of `view' adopted in~\cite{PK11}, as well as here,
  is relevant only to two-dimensional images. 
 By definition, the viewing \textit{direction} is given by
 the vector $w_n\in\BS^2$, which is unique up to sign,
in the representation $HU_n = I-w_nw_n^{\sf T}$, where~$H(x_1,x_2,x_3)=(x_1,x_2,0)$. 
(Thus we regard viewing directions as axial.) 
The `raw view' of the particle~$\rho$ (or its landmarks) is the image of $\rho$ (or its landmarks)
under the map~$HU_n$. We may choose to omit the vanishing third ($x_3$) coordinate, leaving
us with what is called a `profile' in~\cite{Pan09,PK11}. 
The {\sl view\/} of the projected particle (or its landmarks) is obtained 
when the `raw view' (in~$\RR^3$) is identified with (i)~the image of
  any rotation that leaves the projection plane of $H$  
  invariant;~(ii)~the image obtained from reversal of the viewing direction along the $x_3$
  axis. Thus any two images of the particle $\rho$ comprise identical views
   if $\rho$ is subjected to the action $U\rho$ defined as immediately 
   after~\eqref{Pan-model} but now with 
   $U\in\Othree$, rather than~$\SOthree$; and views are equivalence classes 
   with factor group~$\Otwo$, rather than~$\SOtwo$,
   because of view-reversal. See~\cite[Fig.~10, p.~2593]{PK11} for an example.
   

As shown in~\cite[Lemma~5.1, p.~2588]{PK11}, three orthogonal views suffice
to reconstruct the Gram matrix. While, 
as noted in~\cite[pp.~2588--2589]{PK11}, 
it is impossible to ensure that the views selected from the available imagery
are indeed orthogonal,  the procedure they developed fared well enough in their practical example. 
 
 In the fake case, while orientation is no longer Haar distributed, the viewing directions generated by the 
 rotation~$R$ are still uniform. 
 Hence fake uniformity cannot be 
 identified within the tomographic model---unless the modal rotation $M$ from Theorem~\ref{main-thm}
 is  different enough from $I$ to disturb a set of three approximately orthogonal views from the Haar case. 
 Hence fake uniformity does indeed constitute a problem from a modelling viewpoint.

\section*{Acknowledgements}
We are grateful to an anonymous reviewer for helpful comments, which led to 
Section~\ref{ramifications} being completely rewritten, and who pointed out that the distribution of $R$ in Theorem~\ref{main-thm} is 
in general not that of $R^{\sf T}$; see Remark~\ref{law-of-transpose}.

\end{document}